\documentclass[12pt,a4paper]{amsart}
\usepackage[all]{xy}
\usepackage{amssymb}  
\usepackage{latexsym} 

\DeclareFontEncoding{OT2}{}{} 
\newcommand{\textcyr}[1]{%
 {\fontencoding{OT2}\fontfamily{wncyr}\fontseries{m}\fontshape{n}\selectfont #1}}

\newcommand{\Image}{\operatorname{Image}}

\newcommand{\Sha}{{\mbox{\textcyr{Sh}}}}

\newcommand{\Z}{{\mathbb Z}}
\newcommand{\Q}{{\mathbb Q}}

\newcommand{\F}{{\mathbb F}}

\newcommand{\G}{{\mathbb G}}

\newcommand{\To}{\longrightarrow}

\newcommand{\Jac}{\operatorname{Jac}}

\newcommand{\HH}{\operatorname{H}}

\newcommand{\Gal}{\operatorname{Gal}}
\newcommand{\Aut}{\operatorname{Aut}}

\newcommand{\Hom}{\operatorname{Hom}}

\newcommand{\Sel}{\operatorname{Sel}}

\newcommand{\res}{\operatorname{res}}

\newcommand{\Br}{\operatorname{Br}}
\newcommand{\inv}{\operatorname{inv}}

\newenvironment{Proof}{\par\noindent{\sc Proof:}}%
                      {\hspace*{\fill}\nobreak$\Box$\vspace{3mm}\par}
\newenvironment{Remark}{\vspace{3mm}\par\noindent{\sc Remark:}}{\vspace{3mm}\par}%

                      {\hspace*{\fill}\nobreak$\Box$\par}

   {\begin{list}{}{\settowidth{\labelwidth}{#1}%
                   \setlength{\leftmargin}{\labelwidth}%
                   \addtolength{\leftmargin}{\labelsep}}}%
   {\end{list}}

\newtheorem{Theorem}{Theorem}[section]
\newtheorem{Lemma}[Theorem]{Lemma}
\newtheorem{Proposition}[Theorem]{Proposition}
\newtheorem{Corollary}[Theorem]{Corollary}
\newtheorem{Definition}[Theorem]{Definition}
\newtheorem{Question}[Theorem]{Question}

\numberwithin{equation}{section}

\begin{document}
\pagestyle{headings}

\begin{abstract}
Let $A$ be an abelian variety over a number field $k$. We show that weak approximation holds in the Weil-Ch\^atelet group, $\HH^1(k,A)$, but that it may fail when one restricts to the $n$-torsion subgroup. This failure is however relatively mild; we show that weak approximation holds outside a finite set of primes which is generically empty. This proves a conjecture of Land and Tate that can be seen as an analog of the Grunwald-Wang theorem in class field theory. The methods apply, for the most part, to arbitrary finite $G_k$-modules and so may be of interest in their own right.
\end{abstract}

\title{A Grunwald-Wang type theorem for abelian varieties}
\author{Brendan Creutz}
\address{School of Mathematics and Statistics
University of Sydney, Australia}
\email{brendan.creutz@sydney.edu.au}
\maketitle

\section{Introduction}

\subsection{Motivation}Let $k$ be a number field and denote its completion at a prime $v$ by $k_v$. The Grunwald-Wang theorem is an  existence theorem for abelian extensions of $k$ with prescribed local behavior. Namely, given abelian extensions $K_v/k_v$, for $v$ in some finite set $S$, with Galois groups $H_v$, all of which may be embedded in some abelian group $H$, the theorem asserts that there is an abelian extension $K/k$ with completions $K_v$. Usually one is also allowed the requirement that $\Gal(K/k)$ can be embedded in $H$. It is only in a particular set of well understood circumstances that this stronger requirement can fail. For example, the unramified extension of $\Q_2$ of degree $8$ cannot be realized by any degree $8$ cyclic extension of $\Q$. In general, obstructions can only occur when $8$ divides the exponent of $H$, and then only at primes of $k$ lying above $2$.

These global (resp. local) extensions correspond to continuous homomorphisms from the absolute Galois group, $G_k$  (resp. decomposition groups), to $H$. Considering $H$ as a finite $G_k$-module with trivial action, we have $\HH^1(k,H)=\Hom_{cont}(G_k,H)$ and similarly for the local cohomology groups. So the Grunwald-Wang theorem can be rephrased as a statement about the surjectivity of the restriction map $\HH^1(k,H) \to \prod_{v\in S}\HH^1(k_v,H)$.

Now let $A$ be an abelian variety over a number field $k$ and $n$ a positive integer. The Weil-Ch\^atelet group, $\HH^1(k,A)$, parameterizes torsors under $A/k$. In  \cite[p. 683]{LangTate} Lang and Tate write 
\begin{center}
{\em ``In analogy with Grunwald's theorem in class field theory, one may conjecture that if $k$ is an algebraic number field and $\frak{p}$ a given prime, then given $\alpha_{\frak{p}} \in \HH^1(k_\frak{p},A)$, there exists $\alpha \in \HH^1(k,A)$ restricting to $\alpha_{\frak{p}}$.''}
\end{center} For any $n$, there is a surjective map $\HH^1(k,A[n]) \to \HH^1(k,A)[n]$ (here $H[n]$ denotes the $n$-torsion in an abelian group $H$). One can also ask for analogs of the Grunwald-Wang theorem at finite level. 
\begin{Question}
\label{Question1}
Given an abelian variety $A$ over a number field $k$, an integer $n \ge 2$ and a finite set of primes $S$, is the map $\HH^1(k,A[n]) \to \prod_{v\in S}\HH^1(k_v,A[n])$ surjective?
\end{Question}

\begin{Question}
\label{Question2}
Given an abelian variety $A$ over a number field $k$, an integer $n \ge 2$ and a finite set of primes $S$, is the map $\HH^1(k,A)[n] \to \prod_{v \in S}\HH^1(k_v,A)[n]$ surjective?
\end{Question}

Perhaps motivating their conjecture, Lang and Tate showed that the answer to these questions is yes when $\mu_n \subset k$ and the action of the Galois group on $A[n]$ is trivial (\cite{LangTate}, see also \cite[Exercise 10.8]{Silverman}). This also follows directly from the Grunwald-Wang theorem. Assuming the finiteness of the Tate-Shafarevich group $\Sha(A/k)$, Tate went on to characterize the image of the map $\HH^1(k,A) \to \bigoplus_{\text{all }v} \HH^1(k_v,A)$ in terms of a kind of reciprocity coming from the dual abelian variety (see \cite[I.6.26b]{ADT}). However, this settles neither the conjecture nor the questions above, even under the assumption that $\Sha(A/k)$ is finite.

\subsection{Statement of results}
In this paper we prove the conjecture, independently of the finiteness of $\Sha(A/k)$. With regard to questions \ref{Question1} and \ref{Question2} we show, in analogy with the Grunwald-Wang theorem, that the answer can be no in general, but generically will be yes.

\begin{Theorem}[weak approximation]
\label{WCweakapproximation}
Let $A/k$ be an abelian variety over a number field $k$ and let $S$ be any finite set of primes. Then the map $\HH^1(k,A) \to \prod_{v \in S}\HH^1(k_v,A)$ is surjective.
\end{Theorem}

\begin{Theorem}
\label{counterexamples}
There exists an abelian variety $A/\Q$ such that, for infinitely many $n$ (including $n=2$), the map $\HH^1(\Q,A)[n] \to \HH^1(\Q_2,A)[n]$ is not surjective.
\end{Theorem}

\begin{Theorem}
\label{genericwa}
Let $A/k$ be an abelian variety over a number field $k$. There exists a constant $c = c(A,k)$ such that if $n$ is an integer divisible by no prime less than $c$ and $S$ is any finite set of primes, then the map $\HH^1(k,A[n]) \to \prod_{v \in S}\HH^1(k_v,A[n])$ is surjective.
\end{Theorem}

The Grunwald-Wang theorem shows that while weak approximation does not always hold for $\HH^1(k,\Z/n\Z)$, it does hold outside some finite set of primes. Along these lines we prove the following.

\begin{Theorem}[weak weak approximation]
\label{weakweakapproximation}
Let $A/k$ be an abelian variety over a number field $k$ and $n$ an integer. Let $S$ be any set of primes containing all primes of bad reduction and all primes dividing $n$. Let $\Sha(k,A[n],S^c)$ denote the subgroup of $\HH^1(k,A[n])$ consisting of classes that are locally trivial on $S$. The restriction map $\Sha(k,A[n],S^c) \to \prod_{v \notin S}\HH^1(k_v,A[n])$ has dense image in the product of the discrete topologies.
\end{Theorem}

A slightly weaker form of this result (with $\Sha(k,A[n],S^c)$ replaced by $\HH^1(k,A[n])$) can be deduced rather easily from \cite[Lemma I.9.8]{ADT}. In \cite{Creutz} this stronger form is used to show that the $p$-torsion in the Tate-Shafarevich group of any principally polarized abelian variety over a number field is unbounded as one ranges extensions of degree $\mathcal{O}(p)$, the implied constant depending only on the dimension of the abelian variety.

As we show below, these results are closely related to various local-global properties of the group $\HH^1(k,A[n])$. More generally, Dvornicich and Zannier have studied the local-global principle in $G(k)/nG(k)$ and $\HH^1(k,G[n])$, for a commutative algebraic group $G$ (see \cite{DZ1,DZ2,DZ3}). Their results and ideas play a large role in the development below. On the other hand, it seems the related problem of weak approximation has only been addressed in the case of linear algebraic groups (for example \cite[Theorem I.9.10]{ADT} and the Grunwald-Wang theorem itself).

\subsection{Organization}
Section \ref{Preliminaries} contains a quick review of Potiou-Tate duality and other results in Galois cohomology of number fields needed in what follows. In section \ref{LocalGlobal} we study various local-global `principles' for cohomology groups of finite $G_k$-modules. Much of this section is influenced by the ideas in \cite{DZ1,DZ3}. In section \ref{WeakApproximation} we use Poitou-Tate duality to characterize (weak) weak approximation in terms of the local-global principles of the previous section. The main result here is proposition \ref{GW}, which gives a broad generalization of theorem \ref{weakweakapproximation}. All of this is then applied to the particular case of abelian varieties in section \ref{AbelianVarieties}, where the proofs of theorems \ref{WCweakapproximation}, \ref{counterexamples} and \ref{genericwa} are given.\\

\section{Preliminaries}
\label{Preliminaries}

We recall several well known results in Galois cohomology that will be used below. For details we refer the reader to \cite[Section II.5-6]{GC} or \cite[Chapters VII-VIII]{CON}.

Throughout the paper we adopt the following notation:
$k$ is a number field, $G_k$ denotes its absolute Galois group, $M$ is a finite $G_k$-module of exponent $n$ and $M^\vee = \Hom(M,\mu_n)$ is its dual. We use $\HH^i(k,-)$ and $\HH^i(k_v,-)$ to denote global and local Galois cohomology groups. Our convention for indexing (co)products will be that if no index set is specified, then the (co)product is to be taken over all primes $v$ of $k$. Similarly, an expression such as $\prod_{v \notin S}$ is understood to run over all primes of $k$ not in the set $S$. When we say something holds almost everywhere or almost everywhere locally, this means at all but a finite number of primes of $k$.

For each prime $v$ of $k$, duality and the cup-product induce a nondegenerate bilinear pairing (the Tate pairing)
\begin{align}
\label{TatePairing}
(\,, )_v:\HH^1(k_v,M)\times\HH^1(k_v,M^\vee) \to \Br(k_v)[n]\,. 
\end{align} For almost all $v$, the unramified subgroups of $\HH^1(k_v,M)$ and $\HH^1(k_v,M^\vee)$ are exact annihilators with respect to this pairing.

We denote the global pairing induced by the cup product and duality by \[ (\,, ):\HH^1(k,M)\times\HH^1(k,M^\vee) \to \Br(k)[n]\,.\] This pairing is no longer nondegenerate, but it is compatible with the local pairings via the restriction maps. For this reason we will often denote the Tate pairings as defined on global classes as well (i.e. for $\xi \in \HH^1(k,M)$, we may write $(\xi,-)_v$ to mean $(\res_v(\xi),-)_v$ and similarly for classes in $\HH^1(k,M^\vee)$).

For nonarchimedean primes $v$, there is a canonical isomorphism $\inv_v: \Br(k_v) \to \Q/\Z$ (in the archimedean case a canonical injection). The groups $\HH^1(k_v,M)$ and $\HH^1(k_v,M^\vee)$ are finite, and the Tate pairing identifies $\HH^1(k_v,M)$ and $\HH^1(k_v,M^\vee)$ as Pontryagin duals of one another (i.e. $\HH^1(k_v,M^\vee)= \HH^1(k_v,M)^* := \Hom(\HH^1(k_v,M),\Q/\Z)\,$).

The Brauer group of $k$ satisfies a local-global principle expressed by the exactness of
\[ 0 \to \Br(k) \stackrel{\prod \res_v}{\To} \bigoplus \Br(k_v) \stackrel{\sum \inv_v}{\To} \Q/\Z \to 0\,. \]  This gives rise to a product rule for global classes. Namely,
\begin{align}
\label{ProductRule}
\sum \inv_v(\xi,\eta)_v = 0 \text{ for all $\xi\in\HH^1(k,M)$ and $\eta\in\HH^1(k,M^\vee)$.}
\end{align}
Note that since any global class is unramified almost everywhere, this sum is in fact finite.

Given $(\xi_v) \in \prod\HH^1(k_v,M)$, there are two obvious necessary conditions for the existence of a lift of $(\xi_v)$ to a global cocycle. First, $(\xi_v)$ must be unramified outside of some finite set of primes. Second, $(\xi_v)$ must obey the aforementioned product rule. In fact these conditions are already sufficient. This is expressed by the (middle third of the) Poitou-Tate exact sequence
\begin{align}
\label{Poitou1}
\HH^1(k,M) \to \prod '\HH^1(k_v,M) \to \HH^1(k,M^\vee)^*\,. 
\end{align}
Here the product is the restricted product taken with respect to unramified subgroups. The map on the right is given by 
\[ (\xi_v) \mapsto \Bigl( \eta \mapsto \sum\inv_v(\xi_v,\eta) \Bigr) \in \Hom\bigl(\HH^1(k,M^\vee),\Q/\Z\bigr)\,.\]

As a set, the restricted product consists of all families $(\xi_v)$ such that $\xi_v$ is in the unramified subgroup for almost all $v$. It is endowed with a natural topology making it into a locally compact group. The topology is defined by specifying a neighborhood base of $0$ to be the family of all subgroups $\left(\prod_{v \in T} 0\right) \times \left(\prod_{v\notin T} \HH^1_{nr}(k_v,M)\right)$, as $T$ ranges over the finite sets of primes of $k$ containing all primes where $M$ is ramified.

\section{Local-global principles}
\label{LocalGlobal}
\begin{Definition}
\label{Definition1}
Let $M$ be a finite $G_k$-module and $V \subset \HH^1(k,M)$ a subgroup.
\begin{enumerate}
\item For a set $T$ of primes of $k$, we define $\Sha(V,T)$ to be the kernel of the restriction map $V \to \prod_{v \notin T} \HH^1(k_v,M)$.
\item We say that the {\em Hasse principle} holds for $V$ if $\Sha(V,\emptyset)=0$.
\item We say that the {\em strong Hasse principle} holds for $V$ if $\Sha(V,T) = 0$, for every finite set of primes $T$.
\item For a finite set of primes $T$, we say that $V$ is {\em $T$-singular} if the image of the map $\Sha(V,T) \to \prod_{v \in T}\HH^1(k_v,M)$ is not trivial.
\item We say that $V$ is {\em nonsingular} if it is not $T$-singular for any finite set of primes $T$.
\end{enumerate}
\end{Definition}

\begin{Remark}
One easily sees that the strong Hasse principle holds for $V$ if and only if $V$ is nonsingular and the Hasse principle holds. We will see below that the Hasse principle and nonsingularity are, however, independent.
\end{Remark}

If $\xi \in \Sha(V,T)$ for some finite set of primes $T$, we will say that $\xi$ is {\em finitely supported}. If in addition $\xi \notin \Sha(V,\emptyset)$ we say that $\xi$ is {\em $T$-singular} (or simply {\em singular}). When $V = \HH^1(k,M)$ we will use the abbreviation $\Sha(k,M,T)$ for $\Sha(\HH^1(k,M),T)$. For a profinite group $G$ and finite $G$-module $M$, let $\HH^1_*(G,M)$ denote the kernel of the map $\HH^1(G,M) \to \prod_{Z}\HH^1(Z,M)$, where the product runs over all closed cyclic subgroups of $G$. For a Galois extension $K/k$ we will also write $\HH^1_*(K/k,M)$ to denote $\HH^1_*(\Gal(K/k),M)$. This group was introduced by Tate (see \cite{Serre6471} and \cite[Section I.9]{ADT}). This group is also used by by Dvornicich and Zannier to study the Hasse principle for divisibility in commutative algebraic groups, see \cite{DZ1,DZ3}. The following lemma is basically well known (see also \cite[Lemma I.9.3]{ADT}).

\begin{Lemma}
\label{Hcyc}
Let $M$ be a finite $G_k$-module and let $K/k$ denote the minimal Galois extension over which the action on $M$ is trivial. The strong Hasse principle holds for $\HH^1(k,M)$ if and only if $\HH^1_*(K/k,M) = 0$.
\end{Lemma}

\begin{Proof}
Let $T$ be any finite set of primes of $k$ and let $U$ be the set of primes of $K$ which lie above some prime in $T$. The inflation and restriction maps give a commutative diagram with exact rows
\[ \xymatrix{ 0 \ar[r]& \HH^1(K/k,M) \ar[r]^{\inf_{K/k}}\ar[d]^{\prod\res_w}& \HH^1(k,M) \ar[r]^{\res_{K/k}}\ar[d]^{\prod\res_v}& \HH^1(K,M)\ar[d]^{\prod \res_w} \\
0 \ar[r]& \prod_{w \notin U}\HH^1(K_w/k_v,M) \ar[r]& \prod_{w\notin U}\HH^1(k_v,M) \ar[r]& \prod_{w \notin U}\HH^1(K_w,M) }\] (In the bottom row $v$ always denotes the prime of $k$ lying under $w$). The kernel of the middle vertical map is $\Sha(k,M,T)$, while the kernel of the vertical map on the left is $\Sha(\HH^1(K/k,M),T)$ (we are considering $\HH^1(K/k,M)$ as a subgroup of $\HH^1(k,M)$ via the inflation map). The groups on the right consist of continuous homomorphisms. It follows from Chebotar\"ev's density theorem that the vertical map on the right is injective. So the kernels of the other two vertical maps are isomorphic.

It thus suffices to show that $\HH^1_*(K/k,M) = 0$ if and only if the subgroup $\HH^1(K/k,M) \stackrel{\inf_{K/k}}\subset \HH^1(k,M)$ satisfies the strong Hasse principle. This also follows from Chebotar\"ev's density theorem. For all but finitely many primes (namely those where $K/k$ is ramified) the decomposition groups are cyclic and every cyclic subgroup occurs as the decomposition group at a positive density set of primes. So $\HH^1_*(K/k,M)$ consists entirely of finitely supported classes and every finitely supported class in $\HH^1(K/k,M)$ is contained in $\HH^1_*(K/k,M)$. The proof is completed by noting that the strong Hasse principle holds if and only if every finitely supported class is trivial.
\end{Proof}

\begin{Corollary} Let $M$ and $K/k$ be as in the proposition.
\begin{enumerate}
\item The Hasse principle holds for $\HH^1(k,M)$ if and only if $\HH^1_*(K/k,M) \cap \Sha(K/k,M,\emptyset) = 0$.
\item $\HH^1(k,M)$ is nonsingular if and only if $\HH^1_*(K/k,M) \subset \Sha(K/k,M,\emptyset)$.
\end{enumerate}
\end{Corollary}

\begin{Proof} As noted in the proof above, every finitely supported class in $\HH^1(k,M)$ is contained in (the image under the inflation map of) $\HH^1_*(K/k,M)$. The intersection in (1) is trivial if and only if there is no nontrivial class in $\HH^1(k,M)$ with trivial support. Similarly the containment in (2) holds if and only if every finitely supported class has trivial support.
\end{Proof}


From the proof we also extract the following useful observation.
\begin{Corollary}
\label{finitesupport}
Let $M$ be a finite $G_k$-module and $S$ the finite set of primes consisting of all primes where the decomposition group in $\Gal(K/k)$ is not cyclic. Let $V$ be a subgroup of $\HH^1(k,M)$ and $T$ any finite set of primes. Then 
\begin{enumerate}
\item $\Sha(V,T) \subset \Sha(V,S)$.
\item if $T \cap S = \emptyset$, then $\Sha(V,T)=0$.
\end{enumerate}
\end{Corollary}

\begin{Remark}
Note that the set $S$ in the corollary is contained in the finite set of primes where $M$ is ramified.
\end{Remark}

\subsection{Local-global principles for $\HH^1(k,\mu_n)$}
In the case $M = \mu_n$, one can give a complete description of the finitely supported classes in $\HH^1(k,\mu_n)$. Recall that Hilbert's theorem 90 gives an isomorphism $\HH^1(k,\mu_n) \simeq k^\times/k^{\times n}$, so this can be interpreted as the study of local-global properties of divisibility in $\G_m$. We summarize with the following theorem. For the proof we refer the reader to \cite[IX.1]{CON}. The Grunwald-Wang theorem can be derived as a consequence using, for example, theorem \ref{characterize1} below. 

\begin{Theorem}
\label{GWlocalglobal}
Let $T$ be a finite set of primes of $k$, $n = 2^rn'$ be a positive integer with $n'$ odd and let $\kappa$ be the kernel of the map $k^\times/k^{\times n} \to \prod_{v \notin T} k_v^\times/k_v^{\times n}$. Then
\begin{enumerate}
\item $\kappa$ has order dividing $2$.
\item $\kappa$ is nontrivial if and only if $k(\mu_{2^r})/k$ is not cyclic and $T$ contains all primes $v$ which do not decompose in $k(\mu_{2^r})$
\item If $n$ is even, then $\kappa$ is contained in $k^{\times (n/2)}/k^{\times n} \subset k^\times/k^{\times n}$.
\end{enumerate}
\end{Theorem}

\begin{Remark}
Since $k(\mu_4) = k(\sqrt{-1})$ is cyclic, the $\kappa$ can be nontrivial only when $r \ge 3$, i.e. $8 \mid n$. Suppose $k(\mu_{2^r})/k$ is not cyclic and let $S$ be the set of primes of $k$ which do not decompose in $k(\mu_{2^r})/k$. Then $S$ consists entirely of $2$-adic primes. Indeed, all other primes are unramified, so the decomposition groups are cyclic, but $k(\mu_{2^r})|k$ is not.
\end{Remark}

By way of example, consider $V := \HH^1(\Q,\mu_8) \simeq \Q^\times/\Q^{\times 8}$. For a finite set of primes $T$, the kernel of the map $\Q^\times/\Q^{\times 8} \to \prod_{v \notin T }\Q_v^\times/\Q_v^{\times 8}$ is nontrivial if and only if $2 \in T$. Thus the Hasse principle holds for $V$, and the singular sets for $V$ are the finite sets of primes containing $2$. When $2 \in T$, the nontrivial class in the kernel is represented by $16$. In other words $16$ is a $v$-adic $8$-th power if and only if $v \ne 2$. Adjoining a square root of $7$ to $\Q_2$ gives a ramified extension in which $16$ is an $8$-th power. However, $16$ is not an $8$-th power in $\Q(\sqrt{7})^\times$. Thus $\HH^1(\Q(\sqrt{7}),\mu_8)$ is nonsingular but the Hasse principlie fails. Note that, in agreement with theorem \ref{GWlocalglobal}(3), $16$ is a $4$-th power.

\section{Weak approximation}
\label{WeakApproximation}
Suppose $\xi \in \Sha(k,M,T)$ is a class supported entirely on some finite set of primes $T$. Consider its image $\res_T(\xi) \in \prod_{v \in T}\HH^1(k_v,M)$. It follows from the product rule (\ref{ProductRule}) that the image of $\res_T^\vee:\HH^1(k,M^\vee) \to \prod_{v \in T}\HH^1(k_v,M^\vee)$ must be orthogonal to $\res_T(\xi)$ with respect to the nondegenerate pairing \[ (\,, )_T: \prod_{v \in T}\HH^1(k_v,M) \times \prod_{v \in T}\HH^1(k_v,M^\vee) \to \Q/\Z\,,\] given by $\sum_{v\in T}\inv_v(\,, )_v$. In fact, this is the only restriction on the image of $\res_T^\vee$.

\begin{Theorem}
\label{characterize1}
Let $M$ be a finite $G_k$-module with dual $M^\vee$ and $T$ a finite set of primes. An element $(\xi_v) \in \prod_{v\in T}\HH^1(k_v,M^\vee)$ is in the image of $\res_T^\vee$ if and only if it is orthogonal to $\res_T\Sha(k,M,T)$ with respect to the pairing $(\,,\,)_T$. 
\end{Theorem}

The theorem is a special case of proposition \ref{GW}. The more technical version below will allow us to determine also the image of $\res_T^\vee$ modulo arbitrary subgroups of $\prod_{v\in T}\HH^1(k_v,M^\vee)$. This in turn will be used to characterize weak approximation in the $n$-torsion of the Weil-Ch\^atelet group of an abelian variety, and ultimately to prove that weak approximation holds in the Weil-Ch\^atelet group. First we give two corollaries. The second, together with the criterion of Neron-Ogg-Shafarevich, implies theorem \ref{weakweakapproximation}.

\begin{Corollary}
\label{cor1}
The map $\res_T^\vee:\HH^1(k,M^\vee) \to \prod_{v\in T}\HH^1(k_v,M^\vee)$ is surjective if and only if $\HH^1(k,M)$ is not $T$-singular. In particular, weak approximation holds for $\HH^1(k,M^\vee)$ if and only if $\HH^1(k,M)$ is nonsingular.
\end{Corollary}

\begin{Proof}
This follows from the fact that $(\,,\,)_T$ is nondegenerate.
\end{Proof}

\begin{Corollary}
\label{wwa2}
Let $S$ be the finite set of primes consisting of primes where $M$ is ramified.  The map $\res_S^\vee:\Sha(k,M^\vee,S) \to \prod_{v \notin S}\HH^1(k_v,M^\vee)$ has dense image in the product of the discrete topologies.
\end{Corollary}

\begin{Proof}
Suppose $T$ is any finite set of primes and $(\eta_v) \in \prod_{v \in T \cup S}\HH^1(k_v,M^\vee)$ with $\eta_v = 0$, for all $v$ in $S$. We need to show that $(\eta_v)$ is orthogonal to $\res_{T\cup S}\Sha(k,M,T\cup S)$. Clearly $(\eta_v)$ is orthogonal to $\res_{T \cup S}\Sha(k,M,S)$. The result then follows from Corollary \ref{finitesupport} which implies that $\Sha(k,M,T\cup S) \subset \Sha(k,M,S)$.
\end{Proof}

\subsection{Weak approximation for abelian extensions}

\begin{Theorem}[Grunwald-Wang]
Let $T$ be a finite set of primes and, for each $v \in T$, let $K_v/k_v$ be an abelian extension. There exists an abelian extension $K/k$ with completions $K_v$.
\end{Theorem}

\begin{Proof}
Choose an abelian group $A$ for which we can find, for each $v \in T$, an embedding $f_v :\Gal(K_v/k_v) \to A$. It suffices to find an abelian group $B \supset A$ such that $(f_v)$ is in the image of the map $\Hom_{cont}(G_k,B) \to \prod_{v \in T}\Hom_{cont}(G_{k_v},B)$. We reduce to the case that $B$ and, hence, $A$ are cyclic. Suppose $A = \Z/n\Z$ and let $B = \Z/2n\Z$. Let $\xi \in \Sha(k,\mu_{2n},T)$. The dual of $B$ is $\mu_{2n}$, so by theorem \ref{characterize1}, it suffices to show that $(f_v)_{v \in T}$ is orthogonal to $\res_T(\xi) \in \prod_{v \in T}\HH^1(k_v,\mu_{2n})$. By theorem \ref{GWlocalglobal}, $\xi$ lies in the subgroup $n\HH^1(k,\mu_{2n})$. On the other hand, the $f_v$ are $n$-torsion and the pairing $(\,,\,)_T$ is bilinear. The result follows.
\end{Proof}

\begin{Remark}
If every cyclic factor $\Z/n\Z$ of $A$ is such that $\HH^1(k,\mu_n)$ is nonsingular, then one can take $B = A$ in the proof above. Precisely when this can be done is determined by theorem \ref{GWlocalglobal}. Taken together these two results give what is commonly known as the Grunwald-Wang theorem (see \cite[IX.2]{CON}).
\end{Remark}

\subsection{Weak approximation modulo open subgroups}
Let $U$ be an open subgroup of the restricted product $\prod'\HH^1(k_v,M)$. $U$ is a product of subgroups $U_v \subset \HH^1(k_v,M)$. If $U$ is a proper subgroup, then all but finitely many of these are equal to the unramified subgroup. If $U_v^\perp$ denotes the exact annihilator of $U_v$ with respect to the Tate pairing, then $U^\perp := \prod U_v^\perp$ is an open subgroup of $\prod'\HH^1(k,M^\vee)$. Note that when $U = 0$, $U^\perp$ is the entire restricted product and conversely.

Let us use $V_U$ and $V_{U^\perp}$ to denote the subgroups of global classes which map into $U$ and $U^\perp$, repsectively (i.e. \begin{align*}
V_U &= \{ \xi \in \HH^1(k,M)\,|\,\forall v,\, \res_v(\xi) \in U_v\}, \text{ and}\\ V_{U^\perp} &= \{ \eta \in \HH^1(k,M^\vee)\,|\,\forall v,\,\res_v(\eta) \in U_v^\perp\}\,).
\end{align*} 
Let $T$ denote any finite set of primes and let $I' = (\res_T\Sha(V_U,T))^\perp \subset \prod_{v \in T}\HH^1(k_v,M^\vee)$ be the orthogonal complement of $\res_T\Sha(V_U,T)$ with respect to the pairing $(\,,\,)_T$. Let $I$ denote the image of $I'$ under the quotient map $q:\prod_{v \in T}\HH^1(k_v,M^\vee) \to \prod_{v \in T}\frac{\HH^1(k_v,M^\vee)}{U_v^\perp}$.

\begin{Proposition}
\label{GW}
The composition \[ \HH^1(k,M^\vee) \stackrel{\res_T^\vee}{\To} \prod_{v \in T}\HH^1(k_v,M^\vee) \stackrel{q}{\To} \prod_{v \in T}\frac{\HH^1(k_v,M^\vee)}{U_v^\perp}\] maps $\HH^1(k,M^\vee)$ surjectively onto $I$. In particular, $q \circ \res_T^\vee$ is surjective if and only if $V_U$ is not $T$-singular.
\end{Proposition}

Our proof of this proposition is based on (the discussion leading up to) \cite[Theorem 9.2.3]{CON}. Theorem \ref{characterize1} follows by taking $U$ to be the entire restricted product $\prod'\HH^1(k_v,M)$ so that $V_U = \HH^1(k,M)$, $U^\perp = 0$ and $I$ is the orthogonal complement of $\res_T(\Sha(k,M,T))$. We start with a couple lemmas.

\begin{Lemma}
\label{Iproper}
I is a proper subgroup of $\prod_{v \in T}\frac{\HH^1(k_v,M^\vee)}{U_v^\perp}$ if and only if $V_U$ is $T$-singular.
\end{Lemma}

\begin{Proof}
Suppose there is some nonzero element $\xi_T \in \res_T\Sha(V_U,T)$. Since the pairing $(\,,\,)_T$ is nondegenerate, this will be the case if and only if there exists some $\eta_T \in \prod_{v \in T}\HH^1(k_v,M^\vee)$ pairing nontrivially with $\xi_T$. This means $\eta_T \notin I'$. Since $\prod_{v\in T}U_v^\perp$ pairs trivially with $\res_T\Sha(V_U,T)$, this is equivalent to requiring that the class of $\eta_T$ modulo $\prod_{v\in T}U_v^\perp$ does not lie in $I$.
\end{Proof}

\begin{Lemma}
\label{computeIstar}
The Pontryagin dual of $I$ is canonically isomorphic to $\frac{\prod_{v \in T}U_v}{\res_T\Sha(V_U,T)}$.
\end{Lemma}

\begin{Proof}
For finite abelian groups $A_1 \subset A_2 \subset A_3$, with character groups $A_i^*$, let $A_i^\perp$ denote the orthogonal complement of $A_i$ in $A_3^*$ with respect to the natural pairing. One checks that $(A_2/A_1)$ and $A_1^\perp/A_2^\perp$ are canonically identified as duals. The result follows by applying this with \[A_1 = \res_T\Sha(V_U,T) \subset A_2 = \prod_{v\in T}U_v \subset A_3 = \prod_{v\in T}\HH^1(k_v,M)\,,\] since by definition, $I = \frac{\left(\res_T\Sha(V_U,T)\right)^\perp}{\prod_{v\in T}U_v^\perp}$.
\end{Proof}

{\par\noindent{\sc Proof of Proposition \ref{GW}:}}
The discussion leading up to theorem \ref{characterize1} shows that the image of $\res_T^\vee$ is contained in the set $I'$. As $I$ is the image of $I'$ under $q$, the image of the composition in the proposition is contained in $I$. 

By definition $V_{U^\perp}$ is contained in the kernel of $q \circ \res_T^\vee$. Let $\mathcal{C}$ be the cokernel of $q \circ \res_T^\vee$. Taking Pontryagin duals, we have an exact sequence \[ 0 \to \mathcal{C}^* \to I^* \to \Bigl(\frac{\HH^1(k,M^\vee)}{V_{U^\perp}}\Bigr)^*\,. \] Our goal is to show that $\mathcal{C}^*$ is trivial.

The Poitou-Tate exact sequence (\ref{Poitou1}) expresses the fact that an element of $U$ is obtained by restriction of some global class (necessarily in $V_U$) if and only if it is trivial as an element of $\HH^1(k,M^\vee)^*$. Since $U$ pairs trivially with $V_{U^\perp}$ we have an exact sequence \[ V_U \to U \to \Bigl(\frac{\HH^1(k,M^\vee)}{V_{U^\perp}}\Bigr)^*\,.\] This induces an exact sequence
\[ V_U \to \frac{U}{\res\left(\Sha(V_U,T)\right)} \to \Bigl(\frac{\HH^1(k,M^\vee)}{V_{U^\perp}}\Bigr)^*\,.\]

Fitting all of this together, we have a commutative and exact diagram
\[ \xymatrix{ &&0& \\
\Sha(V_U,T) \ar@{^{(}->}[r]& V_U \ar[r]& \prod_{v \notin T}U_v \ar[u] & \\
\Sha(V_U,T) \ar@{^{(}->}[r]\ar@{=}[u] & V_U \ar[r]\ar@{=}[u] & \frac{U}{\res\left(\Sha(V_U,T)\right)} \ar[r]\ar[u] & \Bigl(\frac{\HH^1(k,M^\vee)}{V_{U^\perp}}\Bigr)^* \\
& & \frac{\prod_{v \in T}U_v}{\res_T\left(\Sha(V_U,T)\right)} \ar[u] \ar@{=}[r] &I^*\ar[u]\\
&& 0 \ar[u]& \mathcal{C}^* \ar[u] \\
&&& 0 \ar[u] \\
}\]

The non-tautological equality here is the identification given by lemma \ref{computeIstar}. It follows by a simple chase in the diagram above that $\mathcal{C}^* \simeq \Sha(V_U,T)/\Sha(V_U,T) = 0$.
 {\hspace*{\fill}\nobreak$\Box$\vspace{3mm}\par}

\section{Application to abelian varieties}
\label{AbelianVarieties}

\subsection{Local-global principles for $\HH^1(k,A[n])$}
Let $A$ and $A^\vee$ be dual abelian varieties over $k$ and $n \ge 2$. The following theorem summarizes some of the results of \cite{DZ1,DZ3} regarding the strong Hasse principle in the group $\HH^1(k,A[n])$. Together with corollary \ref{cor1} this proves theorem \ref{genericwa}.

\begin{Theorem}[Dvornicich-Zannier]
\label{DZ}
The strong Hasse principle holds for $\HH^1(k,A[n])$ if any one of the following hold.
\begin{enumerate}
\item $A$ is an elliptic curve and $n$ is prime.
\item $A$ is an elliptic curve over $\Q$ and $n = p^e$ is any power of a sufficiently large prime (independent of $A$).
\item $n = p^e$ is any power of a sufficiently large prime (depending on $A$ and $k$).
\end{enumerate}
\end{Theorem}

\begin{Proof}
Parts (1) and (2) are shown in \cite{DZ1} and \cite[Theorem 1]{DZ3}, respectively. The argument for (3) is given for elliptic curves in \cite{DZ3} (see also \cite[Remark 2.6]{DZ1}). It is not particularly difficult to deduce the same for arbitrary abelian varieties. This goes as follows.

Let $K/k$ be the minimal Galois extension over which the Galois action on $A[p^e]$ is trivial. {\em Suppose} there exists an element $\sigma \in G := \Gal(K/k)$ which acts on $A[p^e]$ as a homothety in $(\Z/p^e\Z)^\times \subset \Aut(A[p^e])$ which has no nontrivial fixed points. In other words $\sigma \in G$ lies in the center of $G$ and $P \mapsto \sigma(P) - P$ is an automorphism of $A[p^e]$. It follows from Sah's lemma (see \cite[Theorem V.5.1]{LangEC}), that $\HH^1_*(G,A[p^e])\subset\HH^1(G,A[p^e])=0$. Then, by lemma \ref{Hcyc}, the strong Hasse principle holds for $\HH^1(k,A[p^e])$.

The existence of such a $\sigma$ follows from a result of Serre \cite{Serre79} (see also \cite[Corollary 2.1.7]{McQuillan}). There exists a constant $d = d(A,k)$ depending only on $A$ and $k$, such that, for any $n$, all $d$-th powers in $(\Z/n\Z)^\times$ arise as homotheties via the action of $G_k$ on $A[n]$. For any prime $p > d+1$ there exists a nontrivial $d$-th power in $\F_p^\times$. Hence, there exists an element $\sigma \in G_k$ which acts on $A[p^e]$ as multiplication by some integer $m_\sigma \ne 1 \mod p$. This implies that $\sigma$ does not fix any nontrivial element of $A[p^e]$.
\end{Proof}

\subsection{Weak approximation for abelian varieties}
We apply the results of the previous section to characterize weak approximation in $\HH^1(k,A[n])$ and $\HH^1(k,A)[n]$. For primes $v$ of $k$, we make the convention that $\HH^0(k_v,A)$ denotes Tate's modified cohomology group; if $v$ is nonarchimedean then $\HH^0(k_v,A) = A(k_v)$, if $v$ is archimedean $\HH^0(k_v,A)$ is the component group of $A(k_v)$ (i.e. the quotient of $A(k_v)$ by the connected component containing the identity).

\begin{Proposition}
\label{WAforAV}
The map $\HH^1(k,A^\vee[n]) \to \prod_{v \in T}\HH^1(k_v,A^\vee[n])$ is surjective if and only if $\HH^1(k,A[n])$ is not $T$-singular. The map $\HH^1(k,A^\vee)[n] \to \prod_{v \in T}\HH^1(k_v,A^\vee)[n]$ is surjective if and only if the $n$-Selmer group of $A$ is not $T$-singular. In particular, weak approximation holds in $\HH^1(k,A^\vee[n])$ if and only if $\HH^1(k,A[n])$ is nonsingular while weak approximation holds in $\HH^1(k_v,A^\vee)[n]$ if and only if $\Sel^{(n)}(A/k)$ is nonsingular.
\end{Proposition}

\begin{Remark}
This implies that weak approximation holds if any of the conditions in theorem \ref{DZ} are met.
\end{Remark}

\begin{Proof}
The statements about $\HH^1(k,A[n])$ follow from the corollaries to theorem \ref{characterize1}. The statement for $\HH^1(k,A)[n]$ will follow by taking $V_U$ in proposition \ref{GW} to be the $n$-Selmer group of $A$.

For any $n$, one has a Kummer sequence, \[ \HH^0(k_v,A) \stackrel{n}{\to} \HH^0(k_v,A) \stackrel{\delta_v}{\To} \HH^1(k_v,A[n]) \to \HH^1(k,A)[n] \to 0\,, \] and similarly for $A^\vee$. The subgroup $U := \prod \Image(\delta_v) \subset \prod\HH^1(k_v,A[n])$ is known to be an open subgroup of the restricted product (the claim is that the image of $\delta_v$ is equal to the unramified subgroup at almost all primes). By definition, the subgroup $V_U \subset \HH^1(k,A[n])$ of classes restricting into $U$ is the $n$-Selmer group of $A$.

Tate's local duality theorems (e.g. \cite[Corollary 3.4]{ADT}) show that the orthogonal complement of $\Image(\delta_v)$ is equal to the image of $\HH^0(k_v,A^\vee)$ under the connecting homomorphism in the Kummer sequence for the dual abelian variety. So proposition \ref{GW} implies that the diagonal map in the commuative diagram below is surjective if and only the $n$-Selmer group of $A$ is not $T$-singular.

\[ \xymatrix{ \HH^1(k,A^\vee[n]) \ar[r]\ar[d]\ar[dr] & \frac{\HH^1(k,A^\vee)}{\delta(\HH^0(k,A^\vee))} \ar@{=}[r]\ar[d] & \HH^1(k,A^\vee)[n] \ar[r]\ar[d] & 0 \\
\prod_{v\in T}\HH^1(k_v,A^\vee[n]) \ar[r]& \prod_{v\in T}\frac{\HH^1(k_v,A^\vee)}{\delta_v(\HH^0(k_v,A^\vee))} \ar@{=}[r]& \prod_{v\in T} \HH^1(k_v,A^\vee)[n] \ar[r] & 0 \\ } \] Clearly the same is true of the vertical map on the right. This is what we wanted to prove.
\end{Proof}

For any $n$, Tate has defined a bilinear pairing 
\[ \langle\,,\,\rangle:\HH^0(k,A)/n\HH^0(k,A)\times\HH^1(k,A^\vee)[n] \to \Br(k)\,.\] This pairing is compatible with the pairing \[ (\,,\,):\HH^1(k,A[n])\times\HH^1(k,A^\vee[n]) \to \Br(k)\] via the Kummer sequences of $A$ and $A^\vee$. Namely, for any $P \in \HH^0(k,A)/n\HH^0(k,A)$ and $\eta \in \HH^1(k,A^\vee)[n]$, $\langle P,\eta \rangle = (\delta(P),\tilde{\eta})$, where $\delta(P)$ denotes the image of $P$ under the connecting homomorphism and $\tilde{\eta}$ denotes any lift of $\eta$ to $\HH^1(k,A^\vee[n])$. The same is true locally and $\langle\,,\,\rangle_v$ indentifies $\HH^0(k_v,A)/n\HH^0(k_v,A)$ and $\HH^1(k_v,A^\vee)[n]$ as Pontryagin duals.

For any finite set of primes $T$, \[ \langle\,,\,\rangle_T := \sum_{v \in T} \inv_v\langle\,,\,\rangle_v\, \] defines a nondegenerate pairing \[ \langle\,,\,\rangle_T:\prod_{v\in T}\frac{\HH^0(k_v,A)}{n\HH^0(k_v,A)}\times\prod_{v \in T}\HH^1(k_v,A^\vee)[n]\to \Q/\Z\,.\] Since, for any $v$, $\res_v\left(\Sel^{(n)}(A/k)\right)$ is contained in the image of $\HH^0(k_v,A)/n\HH^0(k_v,A)$ under the connecting homomorphism, this also gives a pairing
\[ \langle\,,\,\rangle_T : \Sel^{(n)}(A/k) \times \prod_{v\in T}\HH^1(k_v,A^\vee)[n]) \to \Q/\Z\,.\]

Combining proposition \ref{GW} with the proof above readily yields the following.
\begin{Proposition}
\label{WAWCn}
Let $T$ be a finite set of primes. An element $\eta_T \in \prod_{v\in T}\HH^1(k_v,A^\vee)[n]$ is in the image of the restriction map $\HH^1(k,A^\vee)[n] \to \prod_{v \in T}\HH^1(k_v,A^\vee)[n]$ if and only if 
\[ \langle \xi,\eta_T\rangle_T = 0, \text{ for every } \xi \in \Sha(\Sel^{(n)}(A/k),T)\,.\]\\
\end{Proposition}

\subsection{Counter-examples to weak approximation}
Using proposition \ref{WAforAV} we can give examples where weak approximation fails for $\HH^1(k,A[n])$ and $\HH^1(k,A)[n]$.\\

Dvornicich and Zannier \cite{DZ2} have shown that the point $(1561/12^2,19459/12^3)$ on the curve $E : y^2 = (x+15)(x-5)(x-10)$ is divisible by $4$ in $E(\Q_v)$ if and only if $v \ne 2$. This provided one of the first examples of the failure of the strong Hasse principle for $n$-divisibility on an elliptic curve. Under the connecting homomorphism, the point gives rise to a $\{ 2 \}$-singular class in $\Sel^{(4)}(\Q,E) \subset \HH^1(\Q,E[4])$. Weak approximation fails for both $\HH^1(\Q,E)[4]$ and $4$-coverings of $E$ since, by proposition \ref{WAforAV}, $\HH^1(\Q,E)[4] \to \HH^1(\Q_2,E)[4]$ cannot be surjective.\\

Prior to the example above, Cassels and Flynn \cite[p. 61]{CaFlynn} constructed an abelian surface $A$ over $\Q$ for which the Hasse principle for $2$-divisibility fails. We briefly describe their example. Let $C : Y^2 = P(X)Q(X)R(X)$ where $P, Q, R$ are irreducible polynomials of degree $2$ with coefficients in $\Q$ and constant term equal to $1$, splitting over $\Q(\sqrt{2})$, $\Q(\sqrt{17})$ and $\Q(\sqrt{34})$, respectively. The point $a = (0,1) \in C(\Q)$ gives rise to the point $\mathfrak{a} = \{ a, a \} \in A(\Q)$, where $A$ is the Jacobian of $C$.  It follows from their lemma 6.5.1 that $\mathfrak{a} \notin 2A(\Q)$. However, for any $v\le \infty$, at least one of $P,Q,R$ has a root $r_v$ in $\Q_v$. Then the point $\mathfrak{b} = \{ a , (r_v,0) \} \in A(\Q_v)$ is such that $2\mathfrak{b} = \mathfrak{a}$.

Using the same idea, we can construct an example where weak approximation for $\HH^1(\Q,A)[2]$ fails. Let $C$ be the hyperelliptic curve of genus $2$ given by \[y^2 = -(x^2+1)(x^2+5)(x^2-5)\,.\] We have a point $(0,5) \in C(\Q)$. One can check that this gives a $\Q$-point of infinite order $\mathfrak{a} = \{ (0,5),(0,5)\}$ on the Jacobian $A = \Jac(C)$. One easily checks that at least one of the quadratic polynomials defining $C$ has a root in $\Q_v$ if and only if $v \ne 2$. So as above we see that $\mathfrak{a}$ is divisible by $2$ over $\Q_v$ if and only if $v \ne 2$. So the image of $\mathfrak{a}$ in $\Sel^{(2)}(\Q,A) \subset \HH^1(\Q,A[2])$ is $\{2\}$-singular.\\

The next proposition illustrates how the counterexamples above propagate to higher level.

\begin{Proposition}
Suppose that weak approximation fails for $\HH^1(k,A)[n]$. Then weak approximation fails for $\HH^1(k,A)[mn]$ for infinitely many $m$.
\end{Proposition}

\begin{Proof}
We prove the statement for the dual abelian variety. By assumption there is a singular class $\xi \in \Sel^{(n)}(A/k)$. For any positive integer $m$ and any prime $v$, the map $i:A[n] \to A[mn]$ induces a commutative diagram \[ \xymatrix{ \Sel^{(n)}(A/k) \ar[d]^{\res_v}\ar[r]^{i_*}& \Sel^{(mn)}(A/k) \ar[d]^{\res_v} \\
A(k_v)/nA(k_v) \ar[r]^m & A(k_v)/mnA(k_v) }\]  If $\res_v(\xi) = 0$, then $\res_v(i_*\xi) = 0$. So, since $\xi$ is finitely supported, $i_*(\xi)$ is finitely supported. By assumption there exists $v$ such that $\res_v(\xi) \ne 0$. Choose $Q_v \in A(k_v)$ representing $\res_v(\xi) \in A(k_v)/nA(k_v)$. Then $\res_v(i_*(\xi))$ is represented by $mQ_v$. Now suppose $\res_v(i_*(\xi)) = 0$. Then there exists some $P_v \in A(k_v)$ such that $m(Q_v - nP_v) = 0$. Thus $Q_v - nP_v \in A(k_v)[m]$. The torsion subgroup of $A(k_v)$ is finite, so there are infinitely many $m$ for which this cannot happen. For such $m$ we have that $i_*(\xi)$ is singular. The result follows from proposition \ref{WAforAV}.
\end{Proof}

\begin{Remark}
In a similar fashion, Paladino has shown that a counterexample to the Hasse principle in $A(k)/p^nA(k)$ can lead to counterexamples in $A(k)/p^{n+s}A(k)$ for all $s \ge 0$ \cite{Paladino2}. She has also given examples where the Hasse principle fails for divisibility by $9$ in elliptic curves over $\Q$ \cite{Paladino1}. Similar methods give rise to examples of singular classes and the failure of weak approximation in $\HH^1(k,A)[9]$.
\end{Remark}

\subsection{The conjecture of Lang and Tate}
We now come to the proof of theorem \ref{WCweakapproximation}. Since $\HH^1(k,A)$ is torsion, the B\'ezout identity shows that it will suffice to prove the following theorem.

\begin{Theorem}
\label{pprimary}
For any prime number $p$, weak approximation holds in $\HH^1(k,A)[p^\infty]$.
\end{Theorem}
Here $\HH^1(k,A)[p^\infty]$ denotes the subgroup consisting of elements of $p$-powered order. Since we will be working only with $p$-th powers, we use the notation $S^{(n)}(A/k)$ to denote the $p^n$-Selmer group, $\Sel^{(p^n)}(A/k)$. For any positive integers $m$ and $n$, multiplication by $p^{m}$ induces an exact sequence \[0 \to A[p^m] \stackrel{i_*}\to A[p^{m+n}] \stackrel{p^m}\To A[p^{n}] \to 0\,,\] and consequently a map $p^m_*:S^{(m+n)}(A/k) \to S^{(n)}(A/k)$. Let $S(A/k)$ denote the projective limit of the groups $S^{(n)}(A/k)$ with respect to these maps and use $\phi_n:S(A/k) \to S^{(n)}(A/k)$ for the canonical map. One knows that $S(A/k)$ satisfies the strong Hasse principle \cite[Proposition I.6.22]{ADT}. Using this we deduce the following.

\begin{Lemma}
\label{Selnonsingular}
For every $n > 0$, there exits $m$ such that the image of any finitely supported class in $S^{(m+n)}(A/k)$ under the map to $S^{(n)}(A/k)$ is trivial.
\end{Lemma}

\begin{Proof}
We will show below that, for any nontrivial finitely supported class in $S^{(n)}(A/k)$, there exists some $m_0$ such that, for any $m \ge m_0$, no lift of $\xi$ to $S^{(m+n)}(A/k)$ is finitely supported. Using this we prove the lemma as follows. Take $M_0$ to be the maximum of the $m_0$'s as we range over the finitely many nontrivial finitely supported classes in $S^{(n)}(A/k)$ (recall that the $n$-Selmer group is itself finite). Let $M \ge M_0$ and consider a finitely supported class in $S^{(M+n)}(A/k)$. Its image in $S^{(n)}(A/k)$ cannot be equal to any nontrivial finitely supported class. On the other hand, its image is finitely supported. It follows that its image must be trivial.

To establish the claim above, let $L_{\xi} \subset S(A/k)$ be the set of elements which map to $\xi$ in $S^{(n)}(A/k)$. For each $\zeta \in L_\xi$, define $m(\zeta)$ to be the least positive integer such that $\phi_{m(\zeta)+n}(\zeta)$ is not finitely supported. This is well defined since $S(A/k)$ satisfies the strong Hasse principle. Set $m_0 = \sup_{\zeta \in L_\xi}m(\zeta)$. If this supremum is finite the claim follows. So suppose this is not the case. Then we can find a sequence $\{\zeta_i\}_{i=1}^\infty \subset L_\xi$ such that, for each $i$, $\phi_{n+i}(\zeta_i) \in S^{(n+i)}(A/k)$ is finitely supported. Note that $S(A/k)$ is sequentially compact (being a profinite group, it is compact and first countable). So, replacing with a subsequence if necessary, we may assume that the sequence $\zeta_i$ converges to some $\zeta \in S(A/k)$. Clearly $\phi_n(\zeta) = \xi \ne 0$. Our claim will be established if we can show that $\zeta$ is finitely supported, for this will contradict the strong Hasse principle for $S(A/k)$.

Let $S$ be the finite set of primes consisting of primes of bad reduction for $A$ and primes dividing $p$. By the criterion of Neron-Ogg-Shafarevich $A[p^e]$ is unramified outside $S$. Suppose $\res_v(\zeta) \ne 0$ for some $v \notin S$. Then there exists some $M$ such that $\res_v(\phi_{M+n}(\zeta)) \ne 0$. Since $\zeta$ is the limit of the $\zeta_i$, we can choose $i > M$ such that $\phi_{M+n}(\zeta_i) = \phi_{M+n}(\zeta)$ in $S^{(M+n)}(A/k)$. Now $i > M$ and $\phi_{i+n}(\zeta_i)$ is finitely supported, so $\phi_{M+n}(\zeta_i)$ must be as well. By corollary \ref{finitesupport} it follows that $\phi_{M+n}(\zeta) = \phi_{M+n}(\zeta_i)$ is supported on $S$. This contradiction shows that $\zeta$ is supported on $S$, which is a finite set of primes. This is what we intended to show.
\end{Proof}

{\par\noindent{\sc Proof of Theorem \ref{pprimary}:}}
We will prove the statement for the dual abelian variety. Let $T$ be any finite set of primes and $\eta_T := (\eta_v)_{v \in T} \in \prod_{v\in T}\HH^1(k_v,A^\vee)[p^\infty]$. Let $n$ be a positive integer such that, for all $v \in T$, $\eta_v$ is killed by $p^n$. Choose $m$ as in lemma \ref{Selnonsingular} and let $\xi$ be any finitely supported class in the $p^{m+n}$-Selmer group of $A$. We may consider $\eta_T$ as an element in $\prod_{v \in T}\HH^1(k_v,A^\vee)[p^{n+m}]$. Using proposition \ref{WAWCn}, the theorem will follow if we can show that $\eta_T$ is orthogonal to $\xi$ with respect to the pairing \[ \langle\,,\,\rangle_T : S^{(n+m)}(A/k) \times \prod_{v\in T}\HH^1(k_v,A^\vee)[p^{n+m}] \to \Q/\Z\,.\]

We have the commutative diagram 
\[ \xymatrix{ S^{(m+n)}(A/k) \ar[d]^{\res_v} \ar[rr] && S^{(n)}(A/k) \ar[d]^{\res_v}\\ 
\prod_{v \in T}\frac{\HH^0(k_v,A)}{p^{m+n}\HH^0(k_v,A)} \ar[rr] && \prod_{v \in T}\frac{\HH^0(k_v,A)}{p^n\HH^0(k_v,A)} } \] By assumption, the image of $\xi$ in the upper-right hand corner is trivial, so this is also true of its image in the lower-right corner. From commutativity it follows that there exists some $P_T \in \prod_{v \in T}\HH^0(k_v,A)$ such that $\res_T(\xi) \equiv p^nP_T \mod \prod_{v \in T}p^{m+n}\HH^0(k_v,A)$. The pairing $\langle\,,\,\rangle_T$ is bilinear, so \[ \langle \xi,\eta_T\rangle_T = \langle p^{n}P_T, \eta_T \rangle_T = \langle P_T,p^n\eta_T \rangle_T = \langle P_T, 0\rangle_T = 0\,.\] This is what we wanted to show.
{\hspace*{\fill}\nobreak$\Box$\vspace{3mm}\par}

\begin{Remark}
Corollary \ref{finitesupport} shows that any finitely supported class in $\HH^1(k,A[n])$ is unramified outside the set of primes of bad reduction and the primes above $n$. It follows that there can only be finitely many singular classes. The strong Hasse principle is also valid for the projective limit of the groups $\HH^1(k,A[n])$. The proof of \ref{Selnonsingular} carries over for $\HH^1(k,A[n])$ in place of $\Sel^{(n)}(A/k)$. A similar argument to that in the proof above then shows that weak approximation holds for the direct limit of the groups $\HH^1(k,A[n])$ with respect to the maps induced by the obvious inclusion $A[n] \to A[mn]$.
\end{Remark}

\subsection*{Aknowledgements}
I would like to thank Michael Stoll for comments on a preliminary version of this paper and would like to appologize to Roberto Dvornicich and Umberto Zannier for having failed to properly attribute their results in an earlier draft of this paper.


\begin{thebibliography}{MM}

\bibitem[CF]{CaFlynn}
 {\sc J.W.S Cassels and V. Flynn}:
  Prolegomena to a middlebrow arithmetic of curves of genus $2$, 
 {\em London Mathematical Society Lecture Note Series} {\bf 230}, Cambridge University Press, 1996.

\bibitem[Cr]{Creutz}
 {\sc B. Creutz}: Potential $\Sha$ for abelian varieties,
 {\em J. Number Theory}, {\bf 131} (2011) 2162-2174.


\bibitem[DZ1]{DZ1}
  {\sc R. Dvornicich and U. Zannier}:
  Local-global divisibility of rational points in some commutative algebraic groups, {\em Bull. Soc. Math. France} {\bf 129} (2001) 317-338.

\bibitem[DZ2]{DZ2}
  {\sc ---}
  An analogue for elliptic curves of the Grunwald-Wang example,
  {\em C. R. Acad. Sci. Paris} {\bf I 338} (2004) 47-50.

\bibitem[DZ3]{DZ3}
  {\sc ---}
  On a local-global prinicple for the divisibility of a rational point by a positive integer, {\em Bull. London Math. Soc.} {\bf 39} (2007) 27-34.

\bibitem[La]{LangEC}
  {\sc S. Lang}: Elliptic curves: Diophantine analysis, {\em Springer-Verlag}, Berlin and New York, 1978. 

\bibitem[LT]{LangTate}
{\sc S. Lang and J. Tate}:
 Principal homogeneous spaces over abelian varieties,
{\em Amer. Journal of Math.} {\bf 80} (1958) 659-684.

\bibitem[McQ]{McQuillan}
{\sc M. McQuillen}: Division points on semi-abelian varieties,
 {\em Invent. Math.} {\bf 120} (1995) 143-159.

\bibitem[Mi]{ADT}
  {\sc J.S. Milne}:
  Arithmetic duality theorems.
  {\em Perspectives in Mathematics} {\bf 1}, Academic Press, Boston, 1986.

\bibitem[CoN]{CON}
  {\sc J. Neukirch, A. Schmidt and K. Wingberg}:
  Cohomology of number fields (second edition).
  {\em Grundlehren der math. Wissenschaften} {\bf 323}, Springer-Verlag, Berlin, 2000.

\bibitem[Pa1]{Paladino1}
  {\sc L. Paladino}: Elliptic curves with $\Q(E[3])=\Q(\zeta_3)$ and counterexamples to local-global divisibility by $9$, {\em Jour. de Th\'eorie des Nombres de Bordeaux}, {\bf 22} (2010) 139-160.

\bibitem[Pa2]{Paladino2}
  {\sc ---}
  On counterexamples to local-global divisibility in commutative algebraic groups, {\em Acta Arithmetica}, {\bf 148} (2011) 21-29.

\bibitem[Se1]{Serre6471}
  {\sc J.P. Serre}: Sur les groupes de congruence des vari\'et\'es ab\'eliennes, {\em I. Izv. Akad. Nauk SSSR, Ser. Mat.}, {\bf 28} (1964) 3-20; ibid., {\bf 35} (1971) 731-737.

\bibitem[Se3]{Serre79}
  {\sc ---} Quelques propri\'etes des groupes alg\'ebriques commutatifs, Appendix in Ast\'erisque {\bf 69-70} (1979) 191-202.

\bibitem[Se4]{GC}
  {\sc ---}
  Galois cohomology, {\em Springer Monographs in Mathematics}, Springer-Verlag, 2002.

\bibitem[Si]{Silverman}
  {\sc J.H. Silverman}:
  The arithmetic of elliptic curves.
  {\em Graduate Texts in Mathematics} {\bf 106}, Springer-Verlag, 1986.\\

\end{thebibliography}
\end{document}